\documentclass[a4paper,oneside]{amsart}

\usepackage{amssymb,amscd}
\usepackage[applemac]{inputenc}
\usepackage[french]{babel}
\newcommand{\bb}{\mathbb}

\newcommand{\cc}{\bb C}
\newcommand{\rr}{\bb R}

\newcommand{\z}{\bb Z}
\newcommand{\zz}{\z/2}
\newcommand{\pp}{\bb P}

\newcommand{\calO}{\mathcal{ O}}

\newcommand{\xsig}{(X,\sigma)}
\newcommand{\xr}{{X}(\rr)}

\newcommand{\im}{\operatorname{Im}}

\newcommand{\ox}{\mathcal O_X}
\newcommand{\oy}{\mathcal O_Y}

\renewcommand{\ker}{\operatorname{Ker}}

\newcommand{\chiox}{\chi(\ox)}
\newcommand{\chioy}{\chi(\oy)}

\newcommand{\refbooktitle}[1]{\textit{#1}}
\newcommand{\refjourtitle}[1]{\textit{#1}}
\newcommand{\refpapertitle}{}

\numberwithin{equation}{section}%
\newtheorem{theo}[equation]{Théorème}%
\newtheorem{prop}[equation]{Proposition}
\newtheorem{lem}[equation]{Lemme}
\newtheorem{cor}[equation]{Corollaire}

\newtheorem{defin}[equation]{Définition}
\theoremstyle{remark}
\newtheorem{rem}[equation]{Remarque}

\newenvironment{demo}[1]{\noindent\textit{#1.\ }}{\qed\par}

\begin{document}

\title[Surfaces elliptiques réelles]{Nombres de Betti des surfaces elliptiques r\'eelles}
\author{Mouadh Akriche  et Fr\'ed\'eric Mangolte } 

\address{Mouadh Akriche, Laboratoire de Math\'ematiques,
Universit\'e de Savoie, 73376 Le Bourget du Lac Cedex, France,  Fax: +33 (0)4 79 75 81 42} 

\address{Fr\'ed\'eric Mangolte, Laboratoire de Math\'ematiques,
Universit\'e de Savoie, 73376 Le Bourget du Lac Cedex, France, Phone:
+33 (0)4 79 75 86 60, Fax: +33 (0)4 79 75 81 42}
\email{mangolte@univ-savoie.fr}

\begin{abstract} Nous donnons une nouvelle borne pour le nombre de composantes connexes d'une surface elliptique régulière réelle avec section réelle et nous montrons que cette borne est optimale. De plus, toutes les valeurs possibles des nombres de Betti pour de telles surfaces sont atteintes.
\end{abstract}

\maketitle


\begin{quote}\small
\end{quote}

\section{Introduction}\label{sec:intro}

Dans cet article, une surface elliptique est un morphisme $\pi \colon X\longrightarrow \mathbb{P}^1$ où $X$ est une surface analytique complexe compacte, et tel que, à l'exception d'un nombre fini de points $t \in \mathbb{P}^1$,  la fibre $X_t = \pi^{-1}(t)$ est une courbe non singulière de genre 1. 
Sous l'hypothèse que $\pi$ admet au moins une fibre singulière, la surface $X$ est alors régulière, c'est-à-dire qu'elle vérifie $H^1(X,\ox)=\{0\}$. 

La fibration $\pi$ est relativement minimale si aucune courbe exceptionnelle n'est contenue dans une fibre.

Lorsqu'une surface $X$ est munie d'une structure réelle, c'est-à-dire une involution antiholomorphe $\sigma \colon X \to X$, nous disons que le couple  $(X,\sigma)$ est une surface réelle. 
Dans ce cas, l'ensemble des points fixes $X^{\sigma}$ est noté $\xr$, c'est la partie réelle de $(X,\sigma)$. 
Il existe deux structures réelles sur $\mathbb{P}^1$, nous notons $\operatorname{conj} \colon \mathbb{P}^1 \to \mathbb{P}^1$ la conjugaison complexe.
Lorque $\pi \colon X\longrightarrow \mathbb{P}^1$ est elliptique, nous disons que $(X,\sigma)$ (ou plus brièvement $X$) est une surface elliptique réelle si la fibration $\pi \colon X\longrightarrow \mathbb{P}^1$ vérifie $\pi\circ \sigma = \operatorname{conj}\circ \pi$. 
Si $\pi$ admet une section réelle, on en déduit en particulier que $\xr \ne \emptyset$.

Sur $\cc$, les surfaces elliptiques régulières sans fibres multiples sont classifiées par leur caractéristique d'Euler holomorphe $\chiox$. 
Plus précisement, si deux telles surfaces sont relativement minimales, elles sont déformations l'une de l'autre si et seulement si elles ont même caractéristique d'Euler holomorphe, cf.~\cite{Kas77}.

\medskip
Dans cette note, nous montrons que le nombre de composantes connexes de la partie réelle d'une surface elliptique régulière réelle $\pi \colon X\longrightarrow \mathbb{P}^1$ vérifie l'inégalité~:
\begin{equation}\label{eq:comp}
 \# X(\mathbb{R}) \leq 5\chiox 
\end{equation}
lorsque $\pi$ admet une section réelle (Théorème~\ref{theo:bornecompo}).

\medskip
Le premier nombre de Betti de la partie réelle d'une surface elliptique régulière réelle relativement minimale et sans fibres multiples est soumis à l'inégalité~:
\begin{equation}
 h_1(\xr) \leq 10\chiox \;.
\end{equation}

Cette inégalité se déduit directement du Théorème~\ref{theo:rv} dû à Kharlamov qui affirme qu'une telle surface vérifie l'inégalité de Ragsdale-Viro~(\ref{eq:rv}). 
Nous prouvons la majoration (\ref{eq:comp}) grâce à une certaine symétrie sur l'ensemble des nombres de Betti des surfaces elliptiques avec section  réelle, cf. Théorème~\ref{theo:dual}.

\newpage
Dans l'article \cite{Ma00}, le second auteur a montré que l'inégalité de Ragsdale-Viro était optimale pour toute famille complexe de surfaces elliptiques régulières sans fibres multiples. 
Nous montrons dans la Section~\ref{sec:opti} que l'inégalité (\ref{eq:comp}) est optimale elle aussi.

\medskip
Sous les hypothèses considérées, la caractéristique d'Euler holomorphe vérifie $\chiox \geq  1$. Lorsque $\chiox = 1$ (resp. $\chiox = 2$), la surface $X$ est rationnelle (resp. K3). 
Pour ces deux familles, la classification complète des types topologiques de structures réelles est connue, voir l'ouvrage \cite{DIK00} qui contient une bibliographie détaillée.

Pour $\chiox \geq 3$, peu de résultats sont connus. Une surface elliptique régulière $X$ telle que $\chiox \geq 3$ est une surface elliptique propre et en particulier, la fibration $\pi$ est unique. Si $\sigma$ est une structure réelle sur $X$, on a nécessairement  $\pi\circ \sigma = \operatorname{conj}\circ \pi$. 

En appliquant à une surface elliptique régulière réelle avec section réelle les inégalités et congruences habituelles sur les surfaces algébriques réelles cf. e.g. \cite[II.3.9]{Si89} ou \cite{DIK00}, nous obtenons une majoration
\begin{equation}\label{eq:class}
\# X(\mathbb{R}) \leq 5\chiox +\frac {\chiox -\nu}2
\end{equation}
où $\nu \leq 3$ dépend de la classe modulo 8 de $\chiox$.
Pour $\chiox < 4$,  les deux majorations (\ref{eq:class}) et (\ref{eq:comp}) coïncident. Dès que $\chiox \geq 4$, la majoration (\ref{eq:comp}) est meilleure que (\ref{eq:class}).

\medskip
Une partie des résultats exposés dans cet article sont issus de la thèse du premier auteur. Nous remercions vivement J.~Huisman et I. Itenberg pour leurs nombreux conseils. Nous remercions F.~Bihan et J.~van~Hamel pour leurs encouragements durant la thèse.

\section{Modèle de Weierstrass et réduction au cas à fibres nodales}\label{sec:weierstrass}

\bigskip
Pour faciliter la lecture de ce qui suit, nous rappelons une construction classique du modèle de  Weierstrass d'une fibration elliptique avec section \cite{Kas77b}.
Pour que cette construction, à l'origine sur $\cc$, soit valable sur $\rr$, il suffit que la section soit {\em réelle}.

Dans cette partie, une surface elliptique n'est pas forcément lisse, elle peut admettre au pire des points doubles rationnels.  

\bigskip
Soit $\pi \colon X \longrightarrow \mathbb{P}^1$ une surface elliptique régulière réelle. 
Supposons que $\pi$ admette une section réelle $s \colon \mathbb{P}^1 \to X$ et notons  $S = s(\mathbb{P}^1)$ la courbe image.
Munie de la restriction de $\sigma$, $S$ est une courbe définie sur $\rr$.

Supposons la fibration $\pi$ relativement minimale. D'après la classification de Kodaira des fibres singulières possibles \cite[II]{Ko60}, nous savons alors que  la configuration formée par les composantes irréductibles d'une fibre singulière qui ne rencontrent pas $S$ est celle d'une résolution de point double rationnel. 
Notons $\hat X$ la surface obtenue en contractant ces configurations et $\hat \pi \colon \hat X  \longrightarrow \mathbb{P}^1$ la fibration induite, $\hat X$ est alors une surface elliptique éventuellement singulière. 
La section étant réelle, l'ensemble des configurations contractées est globalement stable par $\sigma$ qui induit donc naturellement une structure réelle sur $\hat X$ que nous noterons encore $\sigma$.

Les fibres de $\hat \pi$ sont donc irréductibles et sont  de ce fait isomorphes à des cubiques planes.
Une surface elliptique avec section dont toutes les fibres sont irréductibles est une {\em surface de Weierstrass} et $\hat X$ est le {\em modèle de Weierstrass} de $X$.
Ce modèle est determiné de façon unique aux automorphismes fibrés près.
Par hypothèse, les fibres réelles de $\hat \pi$ possèdent un point réel et sont donc isomorphes sur $\rr$ à des cubiques planes.

\medskip
Notons $k = \chiox$ et $L = (\pi_*N_{S\vert X})^{-1}$ le fibré dual du fibré normal.
Le fibré $L$ est un fibré en droites de degré $k$ sur $\mathbb{P}^1$ qui est défini sur $\rr$. 

Il existe des sections réelles $p \in H^0(\mathbb{P}^1, L^{\otimes 4})$ et $q \in H^0(\mathbb{P}^1, L^{\otimes 6})$ telle que $\hat X$ soit donnée par l'équation 
\begin{equation}\label{eq:cubic}
y^2z = x^3 +pxz^2 +q
\end{equation}
dans $\pp( L^{\otimes 2} \oplus  L^{\otimes 3} \oplus \calO_{\mathbb{P}^1})$.

Par construction, ces sections satisfont les conditions suivantes :
\begin{enumerate}
\item[(i)] Le discriminant $\Delta = 4p^3 + 27 q^2$ (qui est une section de $L^{\otimes 12}$) n'est pas identiquement nul.
\item[(ii)] Pour tout $a \in \pp^1, \min\{3v_a(p), 2v_a(q)\} < 12$  (où $v_a(g)$ est la multiplicité de la section $g$ au point $a$).
\end{enumerate}

Le couple $(p,q)$ est défini de façon unique à l'action de $\cc^*$ suivante près :
\begin{equation}\label{eq:c*action}
\lambda (p,q) = (\lambda^4 p, \lambda^6 q) \;.
\end{equation}

Réciproquement, considérons un triplet de Weierstrass réel $(L,p,q)$. 
C'est-à-dire que $L\to \pp^1 $ est un fibré en droites défini sur $\rr$  et $p \in H^0(\mathbb{P}^1, L^{\otimes 4})$ et $q \in H^0(\mathbb{P}^1, L^{\otimes 6})$ sont des sections réelles qui satisfont les conditions (i) et (ii) précédentes. 
Alors l'équation (\ref{eq:cubic}) définit une surface elliptique réelle $\hat X \to \pp^1$ avec éventuellement des points doubles rationnels dont la résolution minimale $X$ est  une surface elliptique réelle relativement minimale.
La section $S$ est alors la courbe d'équation $x =z =0, y = 1$.

\medskip

Notons $D$ le disque unité  muni de la conjugaison complexe $\operatorname{conj}$. Soit $f \colon \mathcal{M} \longrightarrow D$ une famille de surfaces complexes compactes \cite[I.10]{BPV}. Cette famille est {\em réelle} si 
$\mathcal{M}$ est munie d'une involution anti-holomorphe $\sigma_\mathcal{M}$ telle que $f \circ \sigma_\mathcal{M} = \operatorname{conj} \circ f$. La famille $f \colon \mathcal{M} \longrightarrow D$ est une famille réelle de surfaces elliptiques s'il existe une application holomorphe $\pi \colon \mathcal{M} \longrightarrow  \mathbb{P}^1\times D$ telle que  $\pi \circ \sigma_\mathcal{M} =  \operatorname{conj} \circ \pi$ et qui vérifie les conditions suivantes :
\begin{itemize}
\item $\pi$ admet une section $s$, $s \circ \pi = \operatorname{id}_{\mathbb{P}^1\times D}$, qui commute avec les structures réelles.
\item $\forall t \in D, \pi_t \colon f^{-1}(t) \to \pp^1$ est une surface elliptique dont la fibration elliptique est induite par $\pi$.
\end{itemize}

\begin{defin} Soient $\pi_X \colon X \to \pp^1$ et $\pi_Y \colon Y \to \pp^1$ des surfaces elliptiques réelles avec section réelle.
Nous dirons que $X$ est une {\em déformation réelle} de $Y$ s'il existe une famille réelle de surfaces elliptiques $f \colon \mathcal{M} \longrightarrow D$ et des réels $t \in D$ et $t' \in D$ tels que $\pi_t \colon f^{-1}(t) \to \pp^1$ est isomorphe à $\pi_X \colon X \to \pp^1$ et  $\pi_{t'}\colon f^{-1}(t') \to \pp^1$ est isomorphe à $\pi_Y \colon Y \to \pp^1$.
\end{defin}


\medskip

La plus simple des fibres singulières pouvant apparaître dans une fibration elliptique est une fibre nodale c'est-à-dire une courbe rationnelle avec un point double ordinaire. Lorsqu'une surface elliptique non singulière est de Weierstrass (c'est-à-dire à fibres irréductibles), les seules fibres singulières possibles sont les fibres nodales et les courbes rationnelles avec un cusp.

Soit $\pi \colon X \longrightarrow \mathbb{P}^1$ une surface elliptique régulière réelle relativement minimale avec section réelle. Notons $k = \chiox$. D'après ce qui précède, la surface $X$ est alors dé\-termi\-née par le choix de deux sections réelles $p \in  H^0(\pp^1, \calO_{\pp^1}(4k))$ et $q \in  H^0(\pp^1, \calO_{\pp^1}(6k))$ satisfaisant les conditions (1) et (2), modulo l'action (\ref{eq:c*action}) de $\cc^*$. 
L'ensemble $V$ des couples réels $(p,q) \in H^0(\pp^1, \calO_{\pp^1}(4k)) \times H^0(\pp^1, \calO_{\pp^1}(6k))$ satisfaisant les conditions (i) et (ii) est un ouvert de Zariski non vide.

Grâce à une analyse au cas par cas des types de singularités possibles, Kas \cite{Kas77b} a montré que  la surface $X$ est à fibres singulières nodales si et seulement si le nombre de zéros distincts deux à deux de $\Delta$ est supérieur ou égal à $12k -1$ . 
De telles surfaces forment donc un sous-ensemble ouvert dense de $V$. Nous obtenons en particulier~:

\begin{theo}\label{theo:nodal}
Soit $\pi \colon X \longrightarrow \mathbb{P}^1$ une surface elliptique régulière r\'eelle
relativement minimale. 
Supposons que $\pi$ admet une section r\'eelle, alors $X$ peut être d\'eform\'ee en restant parmi les surfaces elliptiques r\'eelles avec sections r\'eelles en une surface dont les fibres singuli\`eres r\'eelles sont toutes nodales. 
\end{theo}

\section{Symétrie sur l'ensemble des nombres de Betti}

Une fibre singulière nodale complexe est notée $I_{1}$ par Kodaira \cite{Ko60}. On distingue deux types
réels pour une fibre singulière nodale complexe (toujours sous l'hypothèse d'existence d'une section réelle).

\begin{prop} Soit $\pi \colon X \longrightarrow \mathbb{P}^1$ une surface elliptique réelle avec section r\'eelle.
Une fibre singulière réelle $F$ de type complexe $I_{1}$ admet
\begin{itemize}
\item soit deux tangentes réelles en son point double, $F(\rr)$ est alors connexe et nous notons son type réel par $I_{1}^-$, 
\item soit deux tangentes imaginaires conjuguées, $F(\rr)$ est alors non connexe et son type réel est noté $I_{1}^+$.
\end{itemize}
\end{prop}

Voir \cite [Chap.~VII]{Si89} pour les configurations géométriques de tous les types de fibres
singulières réelles d'une fibration elliptique réelle. 

\begin{lem}\label{lem:car}
Soit $\pi \colon X\longrightarrow \pp^1$ une surface elliptique régulière r\'eelle relativement minimale dont les fibres singuli\`eres r\'eelles sont nodales. Alors 
$$
\chi_{top}(\xr) = [I_1^+]-[I_1^-]
$$ 
où $[I_1^+]$ désigne le nombre des fibres singulières de type réel $I_1^+$ et $[I_1^-]$ est le nombre des fibres singulières de type réel  $I_1^-$.
\end{lem}

\begin{demo}{Preuve}
La caractéristique d'Euler topologique de la partie réelle $X(\mathbb{R})$ d'une surface elliptique réelle $\pi \colon X\longrightarrow \pp^1$ est la somme des caractéristiques d'Euler des parties réelles des fibres singulières réelles.
Soit $F$ une fibre singulière réelle nodale, si $F$ est de type réel $I_{1}^+$, $\chi_{top}(F(\rr)) = 1$, sinon, $\chi_{top}(F(\rr)) = -1$.
\end{demo}

\medskip
La topologie de $\xr$ est facile à déterminer dans le cas où les fibres singulières réelles sont toutes, soit de type r\'eel $I_{1}^+$, soit de type réel $I_{1}^-$, car lorsque $\pi$ admet au moins une fibre singulière réelle, la présence d'une section réelle implique qu'au plus une composante connexe est non sphérique, cf. \ref{lem:or}.

\begin{theo}\label{theo:dual}
Soit $\pi \colon X \longrightarrow \mathbb{P}^1$ une surface elliptique régulière réelle, relativement minimale, avec section  r\'eelle. 
Alors, il existe une surface $\pi' \colon X' \longrightarrow \mathbb{P}^1$ elliptique régulière réelle avec section réelle telle que $\chiox = \chi(\mathcal{O}_{X'})$,
$$ 
h_1(X'(\rr)) = 2 \# X(\rr)
$$
et
$$
\# X'(\rr) =\frac{1}{2} h_1(X(\rr))\;.
$$
\end{theo}

\begin{demo}{Preuve du Théorème~\ref{theo:dual}}
D'après le Théorème~\ref{theo:nodal}, nous pouvons supposer que les fibres singulières réelles $X$ sont nodales sans changer de famille complexe ni de type topologique. Il suffit alors de considérer un triplet de Weierstrass réel $(L,p,q)$ associé à $X$, cf. Section~\ref{sec:weierstrass}, et de considérer la surface $X'$ associée au triplet $(L,p,-q)$.  Les surfaces $X$ et $X'$ sont isomorphes sur $\cc$ et la structure réelle de $X'$ est un {\em twist} de la structure réelle de $X$.

Lorsque $X_c = \pi^{-1}(c)$ est une fibre réelle non singulière, $X_c(\rr)$ et  $X'_c(\rr)$ ont le même nombre de composantes connexes. En effet, $\#X_c(\rr)$ est déterminé par le signe du discriminant $\Delta(c) = 4p^3(c) + 27 q^2(c)$.
Lorsque $X_c$ est une fibre réelle singulière, les types réels de $X_c(\rr)$ et $X'_c(\rr)$ sont échangés, $I_1^+$ devient $I_1^-$ et vice versa, cf. e.g. \cite[Chap.~VII]{Si89}. 

Montrons que les sommes totales des nombres de Betti $h_*(\xr)$ et $h_*(X'(\rr))$ sont égales. La partie réelle de la base de la fibration est un cercle et nous pouvons découper ce cercle en arcs $[a,b] \subset \pp^1(\rr)$ tels que $X_a$ et $X_b$ sont des fibres réelles nodales et $X_c$ est une fibre réelle non singulière pour tout $c \in ]a,b[$. 

Notons $arc^+(X)$ (resp. $arc^-(X)$) le nombre d'arcs tels que pour tout $c \in ]a,b[$, $\#X_c(\rr) = 2$ et $X_a$ et $X_b$ sont de type $I_1^+$ (resp. $I_1^-$). Il est clair que $arc^+(X') = arc^-(X)$ et $arc^-(X') = arc^+(X)$. D'autre part, sachant que le nombre de composantes connexes de la partie réelle d'une fibre réelle non singulière change au voisinage d'une fibre réelle nodale, on peut se convaincre sans difficulté que  $h^*(X) = 2 + 2 (arc^+ + arc^-)$. 

Par ailleurs $\chi_{top}(X'(\rr)) = - \chi_{top}(\xr)$ d'après le Lemme~\ref{lem:car} d'où la conclusion sur les nombres de Betti. \end{demo}

\begin{rem}
La deuxième partie de la preuve de \ref{theo:dual} pourrait aussi se prouver à partir de l'invariant $j$ réel défini par J.~Huisman dans \cite{Hu01}.
\end{rem}

\section{Majoration des nombres de Betti}\label{sec:majo}\label{sec:opti}

Soit $\xsig$ une surface algébrique réelle.
Notons $B_i(X) = \dim  H_i(X , \zz )$ le $i$\ieme\ nombre de Betti modulo 2 de $X$ et par
$$
	h_i(\xr) = \dim  H_i(\xr , \zz)
$$
le $i$\ieme\ nombre de Betti modulo 2 de $\xr$. Soit $h_*(\xr) = \sum h_i(\xr)$ et $B_*(X) = \sum B_i(X)$.

\begin{theo}\label{theo:rv}(Kharlamov)
 Soit $\pi\colon X \longrightarrow \mathbb{P}^1$ une surface
 elliptique régulière réelle  sans fibres multiples, alors
 \begin{equation}\label{eq:rv}
 h_1(\xr) \leq h^{1,1}(X) \;.
\end{equation}
\end{theo}

Ce résultat et l'idée de sa preuve ont été communiqués au second auteur par V.~Kharlamov en 1997. Nous n'en connaissons pas de version publiée. Pour le confort du lecteur, nous en proposons ici une rédaction succincte.

\medskip
\begin{demo}{Preuve de \ref{theo:rv}}
Dans le cas où $\pi$ n'admet pas de fibres multiples on vérifie immé\-diate\-ment, en utilisant par exemple la classification des fibres singulières réelles \cite{Si89}, que  
$$
h_1(\xr) \leq  h_1(\operatorname{Jac}(X)(\rr))
$$
où  $\operatorname{Jac}(X)$ est la fibration jacobienne associée à $X$. Par construction, cette fibration est une surface elliptique réelle avec section réelle telle que $h^{1,1}(\operatorname{Jac}(X)) = h^{1,1}(X)$.
Nous supposerons  donc dorénavant que  $\pi$ admet une section réelle sans perdre en généralité.

L'involution $\sigma$ induit une involution $\sigma_*$ sur $H_2(X , \z)$. Considérons les invariants homologiques suivants :

Le rang du sous-module invariant par $\sigma_*$
$$
l^+ = \operatorname{rg} H_2(X , \z)^{\sigma_*} = \operatorname{rg} \ker(1-\sigma_*)
$$
et la caractéristique de Comessatti
$$
 \lambda = \operatorname{rg} \left( (1 + \sigma_*) H_2(X , \z) \right) = \operatorname{rg} \im(1+\sigma_*)\;.
$$

Comme  la fibration $\pi \colon X \to \pp^1$ admet une section, elle n'admet pas de fibres multiples et les nombres de Betti $B_1(X)$ et $B_3(X)$ sont nuls. D'après le théorème de Krasnov, cf. \cite{Kr83}, la surface $\xsig$ est donc {\em Galois-Maximale}. De ce fait, cf. e.g. \cite[Chap. 1]{Si89},
la caractéristique de Comessatti correspond à 
$$
2\lambda = B_*(X) - h_*(\xr)
$$ 
et le premier nombre de Betti de $\xr$ à
\begin{equation}\label{eq:h1}
h_1(\xr) = B_2(X) - l^+ - \lambda \;.
\end{equation}

Si aucune fibre singulière de $\pi$ n'est réelle, $\xr$ est la réunion de deux tores ou de deux bouteilles de Klein car $\pi$ est munie d'une section réelle. Dans ce cas, l'inégalité \ref{eq:rv} est vérifiée. Si $\pi$ admet au moins une fibre singulière réelle, $\xr$ possède une seule composante connexe non simplement connexe et un nombre fini d'autres composantes qui sont toutes homéomorphes à des sphères. Notons $c^+$ le nombre de composantes sphériques. La somme des nombres de Betti de $\xr$ est alors $h_*(\xr) = 2 + 2 c^ + + h_1(\xr)$ et la caractéristique de Comessatti devient
\begin{equation}\label{eq:lambda}
\lambda = l^+ - 2 c^+ \;.
\end{equation}

\medskip
Grâce à un joli raisonnement basé sur la signature de la forme d'intersection, cf. e.g. \cite[Exposé~XIV]{X} ou \cite[3.1.3]{DK00}, $l^+$ est soumis à la majoration
$$
h^{2,0}(X) + c^+ \leq l^+ \;.
$$

De l'équation \ref{eq:lambda}, nous tirons $h^{2,0}(X) - c^+ \leq \lambda$ et la formule \ref{eq:h1} nous permet de conclure~:
$$
h_1(\xr) \leq B_2(X) - 2 h^{2,0}(X)\;.
$$

\end{demo}

\medskip
Dans \cite{Ma00}, le second auteur a montré que la borne précédente est optimale :
 
 \begin{theo}\label{theo:ma}\cite[Théorème~1.2]{Ma00} Toute surface elliptique régulière sans fibres multiples peut être déformée sur $\cc$ en une surface elliptique $X$ admettant une structure réelle telle que 
 $$
 h_1(\xr) = h^{1,1}(X)\;.
 $$
 
 De plus, $X$ admet une fibration elliptique réelle avec une section réelle.
 \end{theo} 

Le premier nombre de Chern d'une surface elliptique relativement minimale $X$ vérifie $c_1^2(X) = 0$ et la formule de Noether donne
\begin{equation}\label{eq:noether}
\chi_{top}(X) = 12\chiox \;.
\end{equation}

Comme  $\pi \colon X \to \pp^1$ n'admet pas de fibre multiple, les nombres de Betti $B_1(X)$ et $B_3(X)$ sont nuls et
\begin{equation}\label{eq:h11}
h^{1,1}(X) = 10 \chiox\;.
\end{equation}

Le théorème principal de cet article est une contrepartie aux deux théorèmes précédents.

\begin{theo}\label{theo:bornecompo}
 Soit $\pi\colon X \longrightarrow \mathbb{P}^1$ une surface
 elliptique r\'eguli\`ere r\'eelle  avec une section  r\'eelle. Alors
 $$
 \# X(\mathbb{R}) \leq 5\chiox \;.
 $$

De plus, pour tout $k \geq 1$, il existe une surface elliptique régulière réelle $X$ avec section  r\'eelle telle que 
$\chiox = k$ et
$$
 \# X(\mathbb{R}) = 5\chiox \;.
 $$
\end{theo}

\begin{demo}{Preuve} La caractéristique d'Euler holomorphe étant un invariant birationnel, il suffit de vérifier l'inégalité pour une fibration relativement minimale. 
En effet, la contraction d'une courbe exceptionnelle réelle ne change pas le nombre de composantes connexes et la contraction d'une paire de courbes exceptionnelles imaginaires conjuguées disjointes non plus. 
Supposons donc que $\pi\colon X \longrightarrow \mathbb{P}^1$ est relativement minimale. D'après le Théo\-rème~\ref{theo:dual}, il existe une surface $X'$ de même caractéristique holomorphe et de premier nombre de Betti $h_1(X'(\rr))=2\#\xr$. Comme $h^{1,1}$ est invariant par déformation, nous déduisons de \ref{theo:rv} la majoration $2\#\xr  \leq h^{1,1}(X)$. D'où la conclusion par \ref{eq:h11}.

Pour montrer que cette borne est atteinte dans une famille complexe donnée, considérons une surface elliptique réelle $\tilde Y$ telle que $\chi(\mathcal{O}_{\tilde Y}) = k$ et $ h_1(\tilde Y(\rr)) = h^{1,1}(\tilde Y)$. L'existence de $\tilde Y$ est assurée par le Théorème~\ref{theo:ma}. Notons $Y$ un modèle relativement minimal de $\tilde Y$. Notons $\eta$ le nombre d'éclatements de points successifs nécessaires pour réaliser $\tilde Y \to Y$. Comme $Y$ admet une fibration elliptique réelle avec section réelle, la symétrie \ref{theo:dual} nous assure qu'il existe une surface  elliptique r\'eguli\`ere r\'eelle relativement minimale $X$ telle que $ \# \xr = \frac 12 h^{1,1}( X) = 5\chi(\mathcal{O}_{Y})$. Pour obtenir finalement une surface dans la famille complexe de $\tilde Y$, éclatons $\eta$ points sur $X$ et  notons $\tilde X$ le résultat. En appliquant un théorème de Kodaira, cf. \cite{Ko60} ou \cite[Th.~5.4]{Ma00}, $\tilde X$ et $\tilde Y$ sont alors dans la même famille complexe. De là, $\chi(\mathcal{O}_{\tilde X}) = \chi(\mathcal{O}_{Y})$ et $\# \tilde X(\rr) = \# \xr = 5 \chi(\mathcal{O}_{\tilde X})$.
\end{demo}

\begin{rem}
Dans sa thèse, le premier auteur à démontré une version plus générale du Théorème~\ref{theo:dual}. Il a obtenu le même résultat de symétrie sur les nombres de Betti pour une fibration elliptique relativement minimale à courbe de base de genre quelconque sous l'hypothèse supplémentaire que l'invariant fonctionnel $J$ soit non constant. En effet, lorsque le genre de la base est non nul,  l'équivalent du Théorème~\ref{theo:nodal} est faux sans cette hypothèse supplémentaire. La preuve, basée sur un théorème de Moishezon \cite{Mo77}, est technique et nous avons préféré utiliser ici une approche plus directe valable uniquement pour les surfaces elliptiques régulières.
\end{rem}

\section{Types topologiques}

Tous les types topologiques possibles pour la partie réelle d'une surface elliptique régulière réelle avec section réelle ont été construits dans \cite{BM06} grâce à une stratégie mêlant  la méthode du patchwork combinatoire de Viro et une technique d'Orevkov  de construction de courbes trigonales via les dessins d'enfants. 

Dans cette partie, nous présentons sans donner de détails, les types topologiques construits par le premier auteur dans sa thèse antérieurement à \cite{BM06}.

\begin{lem}
Le  premier nombre de Betti de la partie réelle d'une surface elliptique régulière réelle relativement minimale et sans fibres multiples est pair.
\end{lem}

\begin{demo}{Preuve}
D'après la théorie de Smith, pour une variété algébrique réelle $X$, la différence $B_*(X) - h_*(\xr)$ est toujours paire. Par ailleurs, pour une surface projective $X$, $B_*(X)$, $\chi_{top}(X)$ et $B_2(X)$ ont même parité. D'après la formule de Noether, cf.~\ref{eq:noether}, $B_2(X)$ est donc pair pour une surface elliptique régulière réelle relativement minimale et sans fibres multiples. Pour terminer, rappelons que $h_1(\xr)$ est de même parité que $h_*(\xr)$.
\end{demo}

\begin{prop}\label{prop:h1possible}
Soit $k$ un entier positif non nul, pour tout entier {\em pair}
$l$, $ l \leq 10k$, il existe  une surface elliptique régulière réelle sans fibres multiples $Y_{k,l}$ de partie réelle connexe et telle que 
$$
\chi(\mathcal{O}_{Y_{k,l}}) = k,\quad  h^{1,1}(Y_{k,l}) = 10k \quad \textrm{ et } \quad h_1(Y_{k,l}(\rr)) = l \;.
$$

De plus, chacune de ces surfaces admet une fibration ellitique avec section réelle dès que $l\ne 0$.
\end{prop}

Par application du Théorème~\ref{theo:dual}, on obtient immédiatement :

\begin{cor}\label{prop:compopossible}
Soit $k$ un entier positif non nul, pour tout entier  $l$,  
$ l\leq 5k$, il existe une surface elliptique r\'eguli\`ere réelle
 $X_{k,l}$ telle que 
$$
 \chi(\mathcal{O}_{X_{k,l}}) = k \quad \textrm{ et } \quad \#X_{k,l}(\rr) = l \;.
$$

De plus, chacune de ces surfaces admet une fibration ellitique avec section réelle dès que $l\ne 0$.
\end{cor}

On note $S_g$ la surface lisse orientable de genre $g$ et
$V_q$ la surface non-orientable de caract\'eristique d'Euler
$2-q$, $V_q$ est difféomorphe à la somme connexe de $q$ copies de $\pp^2(\rr)$.

\begin{prop}\label{prop:top}
Lorsque $k$ est pair, on a $X_{k,l}(\rr) = (l-1)S_0 \sqcup S_1$ et $Y_{k,l}(\rr) = S_l$. Lorsque $k$ est impair, on a  $X_{k,l}(\rr)  = (l-1)S_0 \sqcup V_2$ et $Y_{k,l}(\rr) = V_{2l}$.
\end{prop}

\begin{demo}{Preuve de \ref{prop:top}}
Soit $X$ l'une des surfaces considérées. Cette surface possède une section réelle et au moins une fibre singulière réelle. Sa partie réelle $\xr$ contient donc exactement une composante connexe non simplement connexe $M$ et un nombre fini de composantes homéomorphes à des sphères. Pour déterminer l'orientabilité de $M$ on utilise le lemme suivant qui est un corollaire de \cite[Prop. 6.1]{Ma00}.
\end{demo}

\begin{lem}\label{lem:or} 
Soit $X\to\pp^1$ une surface elliptique régulière réelle relativement minimale avec section réelle, alors $\xr$ est orientable si et seulement si $\chiox$ est pair.
\end{lem}

La Proposition~\ref{prop:h1possible} à été démontrée par le premier auteur dans sa thèse. Pour cela, il a introduit une méthode de transformation de
surfaces elliptiques réelles appellée {\em $I_0^*$-Transformation}. Cette transformation augmente la caractéristique l'Euler holomorphe de 1 tout en permettant un contrôle de  la topologie de la partie réelle.

Cette méthode permet par exemple, de construire une surface elliptique réelle birationnelle à une surface $K3$ à partir d'une surface elliptique réelle rationnelle.

C'est en itérant ce procédé que les surfaces elliptiques réelles de caractéristique d'Euler holomorphe quelconque et à nombres de Betti fixés ont été obtenues à partir de surfaces rationnelles elliptiques réelles.
\medskip

Soit $\pi\colon X\longrightarrow \mathbb{P}^1$ une surface elliptique
 réelle   relativement minimale avec section réelle et $(L,p,q)$ un triplet de Weierstrass réel de $X$. Notons encore $p$ et $q$ les fonctions rationnelles induites par une trivialisation de $L$ au dessus d'un ouvert affine $U$ de $\pp^1$ tel que $\sigma(U)=U$ et au dessus duquel $\pi$ n'admet que des fibres non singulières. Notons $\delta_X = 4p^3+27q^2$ le discriminant. 
 
 Soient $a$ et $b$ deux nombres réels qui ne sont ni des pôles ni des
zéros de $p$ et $q$.
Notons $\pi_Y \colon Y\longrightarrow \mathbb{P}^1$ l'unique modèle non singulier relativement minimal de la surface définie birationnellement par l'équation 
$$
y^2z=x^3+p(u)
\frac{(u-a)^2}{(u-b)^2}xz^2+q(u)
\frac{(u-a)^3}{(u-b)^3}z^3
$$
dans
$\mathbb{P}^2\times U$ munie de la fibration évidente.
La surface obtenue $\pi_Y \colon Y\longrightarrow
\mathbb{P}^1$ est elliptique réelle avec section réelle.

\begin{defin}
La surface $\pi_Y \colon Y\longrightarrow \mathbb{P}^1$ est la
 { \em $I_0^*$-transformée} de $\pi \colon X\longrightarrow \mathbb{P}^1$ par
  $(a,b)$.

\end{defin}

\begin{prop}\label{prop:typereel}
La  $I_0^*$-transformée de $X$ par $(a,b)$ admet
la même liste de fibres singulières complexes que $\pi \colon X \longrightarrow
\mathbb{P}^1$ à laquelle on ajoute exactement 2 fibres singulières
de type $I_0^*$. En particulier, $\chioy = \chiox + 1$.

Soit $\pi^{-1}(c)$ une fibre réelle nodale de $\pi$, avec $c \notin \{a,b\}$, alors $\pi_Y^{-1}(c)$ est de type réel différent de $\pi^{-1}(c)$ si et seulement si $c\in ]a,b[$.
\end{prop}

\begin{demo}{Preuve}
Le discriminant de $Y$ est
$$
\delta_Y= \frac{(u-a)^6}{(u-b)^6}\delta_X
$$ 
et son invariant fonctionnel $j_Y = j_X = 4p^3/27q^2$.

Chaque fibre singulière de $X$ admet le même type topologique qu'une
fibre singulière de $Y$. Les deux fibres singulières de $Y$ au-desssus des
points $u = a$ et $u = b$ sont de type complexe $I_0^*$. Or
$\chi_{top}(I_0^*) = 6$, d'où
$$
\chi_{top}(Y) = 12+\chi_{top}(X) \;.
$$
Nous concluons grâce à la formule de Noether, 
$
\chioy = \chiox + 1\;.
$

\medskip
Notons $F=\pi^{-1}(c)$ et $F'=\pi_Y^{-1}(c)$ les fibres nodales.
Les valuations des fonctions rationnelles $u \mapsto q(u)$ et $u \mapsto q(u).\frac{(u-a)^3}{(u-b)^3}$ en
$c$ sont égales. D'après le Théorème \cite[Theorem~VII.1.5]{Si89}, les types complexes de
$F$ et de $F'$ sont donc les mêmes.
Par ailleurs, les signes de $q(c)$ et  de $q(c).\frac{(c-a)^3}{(c-b)^3}$ sont
opposés si et seulement si $c$ appartient à $[a,b]$. Le type réel de $F'$ est donc 
différent du type réel de $F$ si et seulement si $c\in [a,b]$.
\end{demo}

\medskip

\begin{rem}
Dans \cite{Mi90}, R. Miranda a introduit la notion de $*$-transformée
d'une surface rationnelle elliptique complexe. 
Cette transformation
échange, par exemple, une fibre singulière de type $I_m$ en une fibre
singulière de type $I_m^*$. 
Une $*$-transforma\-tion préserve la
caractéristique d'Euler holomorphe et est donc
différente d'une $I_0^*$-transformation.
\end{rem}


\end{document}